\documentclass[12pt]{article}
\usepackage{amssymb}
\usepackage{latexsym}
\newcounter{myownsection}
\def\myownsection{\refstepcounter{myownsection} \setcounter{equation}{0}}

\begin{document}
$\;$\\[10pt]
\begin{center}
{\bf REDUCTION OF FREE INDEPENDENCE\\
 TO TENSOR INDEPENDENCE}
\footnote{This work is supported by KBN grant No 2P03A00723 and by the EU Network QP-Applications,
Contract No. HPRN-CT-2002-00279}\\[60pt]
{\sc Romuald Lenczewski}\\[30pt]
{\it Institute of Mathematics\\
Wroc{\l}aw University of Technology\\
Wybrze\.{z}e Wyspia\'{n}skiego 27\\
50-370 Wroc{\l}aw, Poland\\
e-mail lenczew@im.pwr.wroc.pl}\\[60pt]
\end{center}
\begin{abstract}
In the hierarchy of freeness construction, free independence was reduced
to tensor independence in the weak sense of convergence of moments.
In this paper we show how to reduce free independence to tensor independence
in the strong sense. We construct a suitable unital *-algebra of closed operators `affiliated' with
a given unital *-algebra and call the associated closure `monotone'. Then we prove that monotone closed operators
of the form
$$
X'= \sum_{k=1}^{\infty}X(k)\overline{\otimes} p_{k}, \;\; X''=\sum_{k=1}^{\infty} p_{k}\overline{\otimes}X(k)
$$
are free with respect to a tensor product state, where $X(k)$ are tensor independent copies of a random variable $X$
and $(p_{k})$ is a sequence of orthogonal projections.
For unital free *-algebras, we construct a monotone closed analog of a unital *-bialgebra called a
`monotone closed quantum semigroup' which implements the additive free convolution,
without using the concept of dual groups.
\\[5pt]
\end{abstract}
\newpage
\myownsection
\begin{center}
{\sc 1. Introduction}
\end{center}
We have shown in [L1] that free independence can be reduced to tensor independence in the following sense.
For a given family of quantum probablity spaces $({\cal A}_{l}, \mu_{l})_{L\in L}$,
there exists a sequence of quantum probability spaces $({\cal A}^{(m)}, \mu^{(m)})_{n\in {\mathbb N}}$
called the {\it hierarchy of freeness} and a sequence of (non-unital) *-homomorphisms
$$
j^{(m)}:\; \sqcup _{l\in L}{\cal A}_{l} \rightarrow {\cal A}^{(m)},
$$
where $\sqcup _{l\in L}{\cal A}_{l}$ is the free product without identification of units, such that we have convergence of moments
$$
\Phi^{(m)}\circ j^{(m)} (X_{1}X_{2}\ldots X_{n})\rightarrow *_{l\in L}\mu_{l}(X_{1}X_{2}\ldots X_{n})
$$
as $m\rightarrow \infty$, where
$X_{1}\in {\cal A}_{l(1)}, \ldots , X_{n}\in {\cal A}_{l(n)}$ and $l(1)\neq l(2)\neq \ldots \neq l(n)$
with $*_{l\in L}\mu_{l}$ denoting the free product of states $\mu_{l}$ in the sense of Avitzour [Av] and Voiculescu [V1].
Moroever, $({\cal A}^{(m)}, \mu^{(m)})_{{m\in \mathbb N}}$ are restrictions
of tensor products of unital *-algebras and states on these *-algebras.

The first order approximation corresponding to $m=1$ gives the boolean product of states [B].
For simplicity, consider two unital *-algebras ${\cal A}$ and ${\cal B}$
and extend them freely by projections $p,p'$ to get $\widetilde{\cal A}={\cal A}*{\mathbb C}[p]$
and $\widetilde{\cal B}={\cal B}*{\mathbb C}[p']$. Then
simple tensors of the form
\begin{equation}
j^{(1)}(X)=X\otimes p' , \;\;\; j^{(1)}(Y) = p\otimes Y
\end{equation}
where $X\in {\cal A}$ and $Y\in {\cal B}$, are boolean
independent with respect to the tensor product state
$\mu^{(1)}=\widetilde{\mu}\otimes \widetilde{\nu}$ on
$\widetilde{\cal A}\otimes \widetilde{\cal B}$, where
$\widetilde{\mu}$, $\widetilde{\nu}$ are the {\it boolean extensions} [L1] of $\mu$ and $\nu$,
respectively.

In the $m$-th order approximation, finite sums of simple tensors of the form
\begin{equation}
j^{(m)}(X)=\sum_{k=1}^{m}X(k)\otimes p_{k}', \;\;\; j^{(m)}(Y)=\sum_{k=1}^{m}p_{k}\otimes Y(k)
\end{equation}
give $m$-free random variables, whose mixed moments of orders $\leq 2m$
in the state $\mu^{(m)}$ agree with moments of free random variables.
Here, $X(k)$'s and  $Y(k)$'s are tensor independent copies of $X$ and $Y$, respectively and
$p_{k}$ as well as $p_{k}'$, $k=1, \ldots , m$, are orthogonal projections.

However, a tensor product representation of free random variables in the strong sense,
i.e. as elements of a tensor product *-algebra, is not so simple since it requires us to take
infinite series of simple tensors instead of finite sums (see [FLS] for the GNS construction).
In this paper we show how this can be done by introducing a suitable closure.
Namely, we adapt to our needs the known concept of the
algebraic closure of a unital *-algebra [Be1-Be3],
which leads us to the notion of the {\it monotone closure}.
As a result, we show that {\it monotone closed operators} of the form
\begin{equation}
j(X)=
\sum_{k=1}^{\infty}X(k)\overline{\otimes}p_{k}', \;\;\;
j(Y)=
\sum_{k=1}^{\infty}p_{k}\overline{\otimes}Y(k)
\end{equation}
are free with respect to a tensor product of states
$\widehat{\mu}$ and $\widehat{\nu}$, which are, roughly speaking, tensor products
of boolean extensions of states $\mu$ and $\nu$ on ${\cal A}$.
Here, $\overline{\otimes}$ denotes the {\it monotone tensor product}, which bears some resemblance
to the von Neumann algebra tensor product.
One can say that in this `quantum orthogonal series'  representation of free random variables,
information about freeness is `shifted from states to variables'.

It is worth pointing out that the representations (1.1)-(1.3) enable us to {\it compare} free random variables
with boolean random variables. In particular, they exhibit in a clear fashion
why units are identified in the free product whereas they are not identified in the boolean product
(nor the $m$-free products).
For finite $m$ we have
$$
\sum_{k=1}^{m}p_{k}\neq 1_{{\cal A}}, \;\;\; \sum_{k=1}^{m}p_{k}'\neq 1_{{\cal B}}
$$
and that is why $j^{(m)}$ does not map the units in ${\cal A}$ and ${\cal B}$
onto the unit in the tensor product and thus we cannot identify them.
In particular, $j^{(1)}(1_{{\cal A}})=1_{{\cal A}}\otimes p'$ and $j^{(1)}(1_{{\cal B}})=p\otimes 1_{{\cal B}}$, respectively.
In turn, the construction of the monotone closure is based on the `completness' property
$$
\sum_{k=1}^{\infty} p_{k}=1_{{\cal A}}, \;\;\; \sum_{k=1}^{\infty}p_{k}'=1_{{\cal B}}
$$
of the sequences of orthogonal projections introduced into the model
and thus $j(1_{{\cal A}})=j(1_{{\cal B}})={\bf 1}$, where ${\bf 1}$ is the unit in the tensor
product. Therefore, in this case units can be identified.

Another important point is to show that the monotone closure can also be introduced on the *-bialgebra level
since coproduct can be viewed as a mapping which produces `independent' copies of a random variable
in a natural and simple fashion.
We already know [L1] how to construct *-bialgebras associated with $m$-freeness if ${\cal A}$ is a unital
free *-algebra generated by a set ${\cal G}$. The simplest example of this type is given by the coproduct
\begin{equation}
\Delta (X)=X\otimes p + p \otimes X
\end{equation}
with $\Delta (p)=p\otimes p$, which produces boolean independent copies of $X\in {\cal G}$ with respect
to the tensor product of extended states, cf. (1.1). Here, $X$
could be called a {\it pre-boolean random variable} when treated as an element of the associated *-bialgebra equipped
with an extended state. In a similar fashion we can produce $m$-free copies of $X$ for all finite $m$.

For $m=\infty$ we introduce the new notion of a {\it monotone closed quantum semigroup},
which is the algebraic structure with *-bialgebra axioms,
in which the algebraic tensor product is replaced by the monotone tensor product.
Now, using the monotone closed quantum semigroup structure on some unital *-algebra
of monotone closed operators ${\cal F}({\cal G})$ associated with ${\cal G}$
we show that by applying the coproduct to monotone closed operators written in the form of a series
\begin{equation}
\sum_{k=1}^{\infty} \delta X(k),
\end{equation}
we obtain, according to (1.3), free copies of $X$, namely
\begin{equation}
\Delta (\sum_{k=1}^{\infty}\delta X(k))= \sum_{k=1}^{\infty}(X(k)\overline{\otimes}p_{k}+
p_{k}\overline{\otimes} X(k))  \;\;({\rm mod}\; {\rm ker}\; \widehat{\mu}\overline{\otimes}\widehat{\nu}),
\end{equation}
which explains why we call the variables of the form (1.5) {\it pre-free random variables}.

This also allows us to reproduce the additive free convolution of states on ${\cal A}$
[V2] using `quantum groups' instead of dual groups [V3]. Recall that the usual convolution of measures
on a group $G$ is implemented by a `quantum group' (Hopf-algebra) structure on some commutative algebra $C(G)$
of functions on $G$. Namely, if $\mu ,\nu$ are functionals on $C(G)$ corresponding to measures on $G$, then
their convolution is given by
$$
\mu \star \nu =(\mu \otimes \nu) \circ \Delta ,
$$
where $\Delta$ is the Hopf-algebra comultiplication on $C(G)$.

For the additive free convolution [V2], an analogous
approach to group-duality was developed by Voiculescu [V3].
Namely, he defined a dual group structure on an algebra ${\cal A}$,
with tensor products replaced by free products. Then
the free convolution of states $\mu$ and $\nu$ on ${\cal A}$ is obtained
from the formula
$$
\mu \boxplus \nu = ( \mu * \nu )\circ \delta
$$
in which the composition of the dual multiplication
$\delta:{\cal A}\rightarrow {\cal A}*{\cal A}$ with the free product
of states $\mu * \nu$
replaces the composition of the Hopf-algebra comultiplication
with the tensor product of states.

Using the monotone closed quantum semigroup structure on ${\cal F}({\cal G})$, we obtain
the additive free convolution $\mu \boxplus \nu$ of states $\mu $, $\nu$ on ${\cal A}$
as a restriction of the quantum semigroup convolution
$$
\widehat{\mu}\star\widehat{\nu}:=(\widehat{\mu}\;\overline{\otimes}
\;\widehat{\nu})\circ \Delta
$$
to the *-subalgebra ${\cal F}_{{\rm pf}}({\cal G})$
of ${\cal F}({\cal G})$ generated by pre-free random variables.
Thus one can view the free additive convolution of classical measures
as a `convolution of quantum measures on a monotone closed quantum semigroup'.

With the results of this paper we complete the program originated in [L1] concerning unification of independence,
or reduction of the main types of independence to tensor independence on extended algebras
(the more recent notion of monotone independence [M] can be also included).
For another general framework, see [L2].

The paper is organized as follows.
In Section 2, we introduce the notion of the monotone closure
for an increasing sequence of unital *-algebras and show that it
has a unital *-algebra structure. In Section 3, we present our main example
of a unital *-algebra ${\cal F}_{0}({\cal G})$ constructed from copies of a unital free *-algebra ${\cal A}$ generated
by the set ${\cal G}$
and we introduce its monotone closure ${\cal F}({\cal G})$. In Section 4 we show that ${\cal F}({\cal G})$ can be endowed
with a monotone closed quantum semigroup structure.
In Section 5, we prove that the associated coproduct produces free random variables which allows us to recover
the free additive convolution from the convolution on the monotone closed quantum semigroup.
In Section 6 we derive a tensor product representation of free random variables.\\[10pt]
\newpage
\myownsection
\begin{center}
{\sc 2. Monotone closed operators}
\end{center}
In this Section we introduce the notion of the monotone closed operators
for certain increasing sequences of unital *-algebras. We are guided by the
construction of the unital *-algebra of closed operators
`affiliated' with a given unital *-algebra [Be1-Be2],
on which we model our notation and terminology.

The construction of the closed *-ring (*-algebra) of operators
consists in taking all sequences
$(x_{m},e_{m})$, where $x_{m}\in {\cal B}_{0}$ and
$(e_{m})$ is a {\it strongly dense domain} (SDD), i.e.
a sequence of projections such that $e_{m}\uparrow 1$ and $x_{n}e_{m}=
x_{m}e_{m}$, $x_{n}^{*}e_{m}=x_{m}^{*}e_{m}$ for $n>m$. On the set of these
sequences called {\it operators with closure} (OWC) one introduces
a suitable equivalence relation, and the set of
the corresponding equivalence classes $[x_{m},e_{m}]$
called {\it closed operators} (CO) `affiliated' with ${\cal B}_{0}$
denoted by ${\cal B}$, can be made into a unital *-ring (*-algebra).
The terminology of this theory is motivated by linear
operators in Hilbert spaces. Heuristically, the ranges of
the $e_{m}$  are an increasing sequence of closed
linear subspaces whose union is a dense linear subspace.
In turn, one can think of $(x_{m},e_{m})$ as a linear
operator whose restriction to the range of $e_{m}$ is $x_{m}e_{m}$.

For instance, this can be done if ${\cal B}_{0}$ is a finite
Rickart *-ring (*-algebra) or, more generally,
a finite Baer *-ring (*-algebra) satisfying LP$\sim$RP,
i.e. the left projection of any $x\in {\cal B}_{0}$ is equivalent
to the right projection of $x$ (equivalence is implemented by a partial
isometry from ${\cal B}$). This procedure can also be applied
to AW$^{*}$-algebras, i.e. Baer *-algebras which are $C^{*}$-algebras
[Be3]. However, a new type of closure is needed for our purposes since
the property LP$\sim$RP is not satisfied in the example
which is of interest to us, namely that
of a unital *-algebra ${\cal F}_{0}({\cal A})$ related to free
products.

Consider an increasing sequence of unital *-algebras
\begin{equation}
{\cal B}^{(0)}\subset {\cal B}^{(1)} \subset {\cal B}^{(2)} \subset \ldots
\end{equation}
where ${\cal B}^{(0)}=\mathbb{C}[p_{1},p_{2}, \ldots]$
is assumed to be the algebra of polynomials
in a countable number of orthogonal projections $(p_{m})$.
Further, we take the union of all algebras ${\cal B}^{(m)}$ denoted
\begin{equation}
{\cal B}_{0}=\bigcup_{m\geq 0}{\cal B}^{(m)}
\end{equation}
and assume that the sequence of increasing projections $(q_{m})$, where
\begin{equation}
q_{m}=p_{1}+p_{2}+\ldots + p_{m},\; m\in {\bf N}
\end{equation}
is an `approximate unit' in ${\cal B}_{0}$, namely
$q_{m}x=xq_{m}=x$ for every $x\in {\cal B}^{(m-1)}$ and $m>1$.

By adding the unit and the zero projection to the sequence $(q_{m})$
we obtain a complete lattice
\begin{equation}
{\cal P}=\{q_{m}; \;0\leq m \leq\infty\}
\end{equation}
where we set $q_{0}=0$ and $q_{\infty}=1$,
which is a sublattice of the lattice of all projections in ${\cal B}^{(0)}$
and thus a sublattice of the lattice of all projections in ${\cal B}_{0}$.
We have $q_{m}\uparrow 1$, both in the sense of lattice supremum in
${\cal P}$ and pointwise in ${\cal B}_{0}$.
Note also that all projections in ${\cal P}$ commute and
the meet of any two projections $e,f\in {\cal P}$ is given by
their product $e\cap f={\rm min}\{e,f\}=ef$.\\
\indent{\par}
{\sc Definition 2.1.}
A {\it monotone strongly dense domain} (MSDD)
in ${\cal B}_{0}$ is a sequence
of projections $(e_{m})$, where $e_{m}\in {\cal P}$
and $e_{m}\uparrow 1$.
If $x\in {\cal B}$ and $e\in {\cal P}$, we write $x^{-1}(e)$
for the largest pojection $g\in {\cal P}$ such that $exg=xg$.
An {\it operator with monotone closure} (OWMC)
is a sequence $(x_{m},e_{m})$
with $x_{m}\in {\cal B}^{(m)}$ and $(e_{m})$ a MSDD, such that $m<n$ implies
$x_{n}e_{m}=x_{m}e_{m}$ and $x_{n}^{*}e_{m}=x_{m}^{*}e_{m}$.\\
\indent{\par}
{\it Example.}
Of course, $(q_{m})$ is a MSDD, but an equally important example
for us will be the shifted sequence $(q_{m-k})=(e_{m})$, where
$$
e_{m}
=
\left\{
\begin{array}{cc}
0 & {\rm if} \; m\leq k\\
q_{m-k} & {\rm if} \; m>k
\end{array}
\right.
$$
i.e. $(q_{m-k})$ is also a MSDD -- in the sequel the notation
$(q_{m-k})$ will be used with the
understanding that the index $m$ is reserved for the running index,
whereas $k$ is fixed and is responsible for the shift.
Similarly, the sequence $(f_{m})=(1_{m-k})$ defined by
$$
f_{m}
=
\left\{
\begin{array}{cc}
0 & {\rm if} \; m\leq k\\
1 & {\rm if} \; m>k
\end{array}
\right.
$$
is also a MSDD. This MSDD is used to embed ${\cal B}_{0}$ in ${\cal B}$.
Namely, if $x\in {\cal B}^{(k)}$, then the MCO $[x_{m},1_{m-k}]$, where
$$
x_{m}=
\left\{
\begin{array}{cc}
0 & {\rm if}\;\;m\leq k\\
x & {\rm if}\;\; m>k
\end{array}
\right.
$$
can be identified with $x$ (all such embeddings are
consistent since they have the same `tails').\\
\indent{\par}
{\it Remark.}
It is sometimes convenient to write OWMC (and MCO) in the form of series
$$
(x_{m},e_{m})=\sum_{m=1}^{\infty}(x_{m}-x_{m-1})
$$
where $x_{m}\in {\cal B}^{(m)}$ and we set $x_{0}=0$
(we then keep in mind the SDD $(e_{m})$).\\
\indent{\par}
{\sc Lemma 2.2.}
{\it Suppose $(e_{m})$ and $(f_{m})$ are} MSDD.
{\it Then $(e_{m}f_{m})$ is an} MSDD. {\it Further, if a sequence $(x_{m})$,
where $x_{m}\in {\cal B}^{(m)}$ for every $m$,
satisfies $x_{n}e_{m}=x_{m}e_{m}$ for every $n>m$,
then the sequence $(g_{m})=(e_{m}x_{m}^{-1}(f_{m}))$
is a} MSDD.\\[5pt]
{\it Proof.}
That $(e_{m}f_{m})$ is a MSDD if $(e_{m})$ and $(f_{m})$ are MSDD, immediately
follows from the definition of MSDD.
Let us prove that $(g_{m})$ is a MSDD.
The proof of monotonicity is the same as in the case
of SDD and we quote it after [Be1-Be2] only for the reader's convenience.
Denote $h_{n}=x_{n}^{-1}(f_{n})$. Thus $g_{n}=e_{n}h_{n}$ and $h_{n}$
is the largest projection from the lattice  ${\cal P}$ such that
$$
(1-f_{n})x_{n}h_{n}=0
$$
If $n>m$, then
$$
x_{n}g_{m}=x_{n}e_{m}g_{m}=x_{m}e_{m}g_{m}=x_{m}g_{m}=x_{m}h_{m}g_{m}
$$
and thus
$$
(1-f_{n})x_{n}g_{m}=(1-f_{n})x_{m}h_{m}g_{m}=0
$$
which, by maximality of $h_{m}$, implies that
$g_{m}\leq h_{n}$. Also, $g_{m}\leq e_{m}\leq e_{n}$ since
$(e_{m})$ is a MSDD. Thus, $g_{m}\leq h_{n}e_{n}=g_{n}$, which ends the proof
of monotonicity of $(g_{m})$.

Let us show that $g_{m}\uparrow 1$. Since $(f_{m})$ is a MSDD, we have
$f_{m}=q_{l(m)}$ where $0\leq l(m) \leq \infty$ and
$l(m)\uparrow \infty$.
We now use the assumption that each $x_{k}\in {\cal B}^{(k)}$,
which implies that $f_{m}$ acts as an identity when
multiplied by $x_{1}, x_{2}, \ldots , x_{l(m)-1}$.
Hence, if $l(m)>m$, then $(1-f_{m})x_{m}=0$, which gives
$h_{m}=1$.
In turn, if $l(m)\leq m$, then we write
$$
x_{m}=x_{l(m)-1}+(x_{m}-x_{l(m)-1})
$$
and we have
\begin{equation}
(1-f_{m})x_{m}h_{m}= (1-f_{m})(x_{m}-x_{l(m)-1})h_{m}
\end{equation}
since
$$
f_{m}x_{m}=x_{l(m)-1}+f_{m}(x_{m}-x_{l(m)-1})
$$
but equation
$$
(1-f_{m})(x_{m}-x_{l(m)-1})h_{m}=0
$$
is satisfied if we take for $h_{m}$ any projection
which right-annihilates $(x_{m}-x_{l(m)-1})$.
Note that $e_{l(m)-1}$ is such a projection by assumption. Therefore,
$h_{m}\geq e_{l(m)-1}$. Since $l(m)\rightarrow \infty$ and
$(e_{m})$ is a MSDD, we obtain $g_{m}=h_{m}e_{m}\uparrow 1$.\\
\indent{\par}
{\it Remark.}
One should point out two new features in our definition of OWMC
as compared to OWC: we assume that
$x_{m}\in {\cal B}^{(m)}$
and we take a  more restricted family of SDD.
Note that in the theory of closed *-rings (*-algebras)
we have $y^{-1}(e)\succeq e$ for any projection $e$ and any $y$,
where $\succeq$ is the order implemented by partial isometries.
Using this and LP$\sim$RP, one shows that
$e_{m}\cap x_{m}^{-1}(f_{m})$ is a SDD under assumptions similar to
those in Lemma 2.2.
However, in the main example studied in this
paper, the property $y^{-1}(e)\succeq e$ does not hold, which is one of the
obstacles in applying the usual theory.\\
\indent{\par}
{\sc Lemma 2.3.}
{\it If  $(x_{m},e_{m})$ and $(y_{m},f_{m})$ are} {\rm OWMC}
{\it and we define
\begin{equation}
k_{m}=f_{m}y_{m}^{-1}(e_{m})e_{m}(x_{m}^{*})^{-1}(f_{m})
\end{equation}
then $(x_{m}^{*},e_{m})$, $(x_{m}+y_{m},e_{m}f_{m})$,
$(\lambda x_{m},e_{m})$ and $(x_{m}y_{m}, k_{m})$ are} {\rm OWMC},
{\it where $\lambda\in \mathbb{C}$.}\\[5pt]
{\it Proof.}
Note that $x_{m}y_{m}\in {\cal B}^{(m)}$ since
$x_{m},y_{m}\in {\cal B}^{(m)}$.
Then one shows that the sequences $(x_{m}+y_{m},e_{m}f_{m})$,
$(x_{m}y_{m},k_{m})$, $(\lambda x_{m}, e_{m})$ and
$(x_{m}^{*},e_{m})$ are OWMC. For instance
\begin{eqnarray*}
x_{n}y_{n}k_{m}&=&
x_{n}y_{n}e_{m}f_{m}y_{m}^{-1}(e_{m})(x_{m}^{*})^{-1}(f_{m})\\
&=& x_{n}y_{m}y_{m}^{-1}(e_{m})e_{m}f_{m}(x_{m}^{*})^{-1}(f_{m})\\
&=& x_{n}e_{m}y_{m}k_{m}\\
&=& x_{m}e_{m}y_{m}k_{m}\\
&=& x_{m}y_{m}k_{m}
\end{eqnarray*}
where we use the definition of $y_{m}^{-1}(e_{m})$
and the fact that all projections involved commute.
\hfill $\Box$\\
\indent{\par}
{\sc Definition 2.4.}
We say that the OWMC $(x_{m},e_{m})$, $(y_{m},f_{m})$
are {\it equivalent}, written $(x_{m},e_{m})\equiv (y_{m},f_{m})$ if there
exists an MSDD $(g_{m})$ such that $x_{m}g_{m}=y_{m}g_{m}$ and
$x_{m}^{*}g_{m}=y_{m}^{*}g_{m}$ for every $m$.
This relation is an equivalence relation and
we say that $(g_{m})$ {\it implements} the equivalence.
The set of all equivalence classes
is denoted by ${\cal B}$ and its elements are called {\it monotone
closed operators} (MCO).
We denote by $[x_{m},e_{m}]$ the MCO determined by the OWMC
$(x_{m},e_{m})$.\\
\indent{\par}
{\sc Theorem 2.5.}
{\it If ${\bf x}=[x_{m},e_{m}], {\bf y}=[y_{m},f_{m}] \in {\cal B}$, then
the operations
\begin{eqnarray*}
{\bf x}+{\bf y}&=&[x_{m}+y_{m},e_{m}f_{m}]\\\;
{\bf x}{\bf y}&=&[x_{m}y_{m},k_{m}]\\\;
{\bf x}^{*}&=&[x_{m}^{*},e_{m}]\\\;
\lambda {\bf x}&=&[\lambda x_{m},e_{m}]
\end{eqnarray*}
are well-defined and make ${\cal B}$ into a unital *-algebra, where
$k_{m}$ is given by (2.2) and $\lambda\in \mathbb{C}$.}\\[5pt]
{\it Proof}.
The proof since it is similar to that in the case of
closed operators [Be1,Be2].\hfill $\Box$\\[10pt]
%
\myownsection
\begin{center}
{\sc 3. The unital *-algebra
${\cal F}_{0}({\cal G})$ and its monotone closure
${\cal F}({\cal G})$}
\end{center}
\nopagebreak[4]
In this Section we give a construction of a unital *-algebra that is
related to free products of states in free probability.
The main idea of the
construction is to extend the free product of a countable number of copies
of a given unital *-algebra by a sequence of projections which will
provide us with a `nice' SDD.

Let ${\cal A}$ be the unital free *-algebra generated by the set ${\cal G}$.
We take countably many copies of
this algebra which we label by ${\cal A}'(k)$ and
${\cal A}''(k)$, where $k\in {\bf N}$, with the condition
${\cal A}''(1)=\{0\}$.
Now, take the free product of
all these copies (without identification of units)
$$
\widehat{\cal F}({\cal G})
=\sqcup _{k\in {\bf N}}({\cal A}'(k)\sqcup {\cal A}''(k))
$$
i.e. the linear span of all words
(we treat the unit from $\mathbb{C}$ as the empty word),
and extend this product by the *-algebra of polynomials over
the sequence $(p_{m})$ of orthogonal projections, namely
$$
\widehat{\cal F}_{0}({\cal G})
=\widehat{\cal F}({\cal G})*\mathbb{C}[p_{1},p_{2}, p_{3}, \ldots ]
$$
where $p_{k}p_{l}=\delta_{k,l}p_{k}$, $p_{k}^{*}=p_{k}$ and we assume that
the unit in $\widehat{\cal F}({\cal G})$ (the empty word)
and the unit in $\mathbb{C}[p_{1},p_{2}, p_{3}, \ldots ]$ are identified and
denoted by $1$.

In other words, $\widehat{\cal F}_{0}({\cal G})$ can be defined as
the linear span of words of the form
\begin{equation}
s_{1}w_{1}s_{2}w_{2}\ldots s_{n}w_{n}s_{n+1}
\end{equation}
with the juxtaposition product,
where $w_{1}, \ldots ,w_{n}$
are non-empty words from  $\widehat{\cal F}({\cal G})$ and
$s_{1}, \ldots , s_{n+1}$ are non-trivial projections from the lattice
${\cal P}$ given by (2.4).

For given $X\in {\cal G}$, we denote by
$X'(k)$ and $X''(k)$ the copies
of $X$ in ${\cal A}'(k)$ and ${\cal A}''(k)$,
respectively.\\
\indent{\par}
{\sc Definition 3.1.}
Let $J$ be the two-sided *-ideal in $\widehat{\cal F}_{0}({\cal G})$
generated by elements of the form
\begin{eqnarray}
(1-q_{m})X'(k)&=&0 \;\; {\rm for}\;\; k<m\\
(1-q_{m})X''(k)&=&0\;\; {\rm for}\;\; k<m\\
q_{m}(X'(k)-X''(k))&=&0\;\; {\rm for}\;\; k>m
\end{eqnarray}
where $X\in {\cal G}$. We denote by
${\cal F}_{0}({\cal G})=\widehat{\cal F}_{0}({\cal G})/{\cal J}$ the corresponding
quotient algebra.\\
\indent{\par}
In other words, we can define ${\cal F}_{0}({\cal G})$ as
the linear span of {\it reduced} words, i.e. words of the form (3.1) with the minimal
number of non-trivial projections from the lattice ${\cal P}$
(after (3.2)-(3.4) have been taken into account) with the
juxtaposition product inherited from ${\cal F}_{0}({\cal G})$.

If we denote by
${\cal F}^{(m)}({\cal G})$ the unital *-subalgebra of
${\cal F}_{0}({\cal G})$ spanned by reduced words built from all projections
$(q_{m})$ and generators from ${\cal A}'(k), {\cal A}''(k)$
with $1\leq k\leq m$ and
${\cal F}^{(0)}({\cal G})=\mathbb{C}[p_{1},p_{2}, \ldots ]$,
then the sequence $({\cal F}^{(m)}({\cal G}))$ is an increasing sequence
of unital *-algebras (2.1) with union ${\cal F}_{0}({\cal G})$ as in (2.2),
for which monotone closure can be constructed
along the lines of Section 2. In particular,
relations (3.2)-(3.3) make $(q_{m})$ into a MSDD in ${\cal F}_{0}({\cal G})$.
In turn, relation (3.4) will be needed when introducing
certain OWMC associated with free random variables.
By ${\cal F}({\cal G})$ we denote the monotone closure of ${\cal F}_{0}({\cal G})$ consisting
of monotone closed operators $[x_{m},e_{m}]$, where $x_{m}\in {\cal F}_{0}({\cal G})$
and $(e_{m})$ is a MSDD in ${\cal F}_{0}({\cal G})$.\\
\indent{\par}
{\sc Definition 3.2.}
The monotone closed operators of the form
\begin{equation}
[x_{m},q_{m}]=[\sum_{k=1}^{m}\delta X(k), q_{m}]
\end{equation}
where
\begin{equation}
\delta X(k)=
\left\{
\begin{array}{cc}
X'(k)-X''(k) &{\rm if} \;\;k>1\\
X'(1)&{\rm if}\;\; k=1
\end{array}
\right.
\end{equation}
and $X\in {\cal G}$, will be called {\it pre-free random variables}.\\
\indent{\par}
Pre-free random variables contain encoded information about freeness. Namely, after applying
an appropriate comultiplication to a pre-free random variable we obtain
a sum of free random variables with respect to a tensor product of states (see Section 5).
Denote by ${\cal F}_{{\rm pf}}({\cal G})$ the (non-unital) *-subalgebra
of ${\cal F}({\cal G})$ generated by pre-free random variables.

To get a glimpse of freeness in this definition, let us now derive the explicit form of the product of
pre-free random variables
(the remaining algebraic operations are straightforward).\\
\indent{\par}
{\sc Theorem 3.3.}
{\it The product of pre-free random variables takes the form}
$$
[x_{m}^{(1)},q_{m}][x_{m}^{(2)},q_{m}]\ldots [x_{m}^{(n)},q_{m}]
=[x_{m}^{(1)}x_{m}^{(2)}\ldots x_{m}^{(n)},q_{m-n+1}]
$$
{\it where $[x_{m}^{(k)},q_{m}]\in
{\cal F}_{{\rm pf}}({\cal G})$, $k=1, \ldots , n$,
are pre-free random variables associated with $X_{1}, \ldots ,
X_{n}\in {\cal G}$.}\\[5pt]
{\it Proof.}
Instead of giving a formal inductive proof, we prefer to analyze
the cases $n=2$ and $n=3$, which is more intuitive and seems sufficient to
see how to proceed in the general case.
Let $X,Y,Z\in {\cal G}$ be non-zero and
let $[x_{m},q_{m}]$, $[y_{m},q_{m}]$, $[z_{m},q_{m}]$ be the associated pre-free random variables.
We begin with showing that
$$
[x_{m},q_{m}][y_{m},q_{m}]=
[x_{m}y_{m}, q_{m-1}].
$$
It can be seen from (2.6) that this boils down to the computation of
$y_{m}^{-1}(q_{m})$.
We have to find the largest projection $g_{m}\in {\cal P}$ such that
\begin{equation}
(1-q_{m})\sum_{k=1}^{m}\delta Y(k)g_{m}=0
\end{equation}
Note that
$$
(1-q_{m})\sum_{k=1}^{m-1}\delta Y(k)=0
$$
for any $m>1$ in view of (3.2)-(3.3).
Hence, (3.7) holds iff
$$
(1-q_{m})\delta Y(m)g_{m}=0
$$
and the largest projection from ${\cal P}$ which satisfies this equation
is the largest projection from ${\cal P}$ which right-annihilates
$\delta Y(m)$, namely $q_{m-1}$. Thus, $g_{m}=q_{m-1}$.
Similarly, $(x_{m}^{*})^{-1}(q_{m})=q_{m-1}$ and
therefore, the $k_{m}$ of (2.6) are given by $k_{m}=q_{m-1}$.
The next step consists in finding the $(k_{m})$ in the formula
$$
(z_{m},q_{m})(x_{m}y_{m},q_{m-1})=(z_{m}x_{m}y_{m},k_{m}),
$$
which is of the form
$$
k_{m}=q_{m-1}(x_{m}y_{m})^{-1}(q_{m})q_{m}(z^{*}_{m})^{-1}(q_{m-1}).
$$
First, let us find the largest $g_{m}\in {\cal P}$ such that
\begin{equation}
(1-q_{m})x_{m}y_{m}g_{m}=0.
\end{equation}
We claim that if we take $g_{m}=q_{m-2}$, then this equation is satisfied.
Namely,
$$
x_{m}y_{m}q_{m-2}=x_{m}y_{m-2}q_{m-2}=
x_{m}q_{m-1}y_{m-2}q_{m-2}=x_{m-1}q_{m-1}y_{m-2}q_{m-2}
$$
where we used (3.2)-(3.4) repeatedly, and thus
$$
q_{m}x_{m}y_{m}q_{m-2}=x_{m-1}q_{m-1}y_{m-2}q_{m-2}=
x_{m}y_{m}q_{m-2}
$$
which  implies that $g\geq q_{m-2}$. We need to show that
$q_{m-2}$ is the largest projection from ${\cal P}$ which
satisfies (3.8). Let us check if $q_{m-1}$ solves (3.8) when
substituted for $g_{m}$. Using (3.2)-(3.4) we get
\begin{eqnarray*}
(1-q_{m})x_{m}y_{m}q_{m-1}&=&(1-q_{m})x_{m}y_{m-1}q_{m-1}\\
&=&
(1-q_{m}) \delta X(m)y_{m-1}q_{m-1}\\
&=&
\sum_{l=1}^{m-1}(1-q_{m})\delta X(m)\delta Y(l)q_{m-1}\neq 0
\end{eqnarray*}
since we obtained a sum of linearly independent words, each of which
is non-zero.
In a similar way we show that $g_{m}\neq q_{m}$ and
it is clear that $g_{m}\neq q_{p}$ for $p>m$ since
$(1-q_{m})x_{m}y_{m}q_{p}=(1-q_{m})x_{m}y_{m}\neq 0$.
This finishes the proof for the product of three OMC. For products
of higher order we continue in a similar way.
\hfill $\Box$\\
\indent{\par}
{\it Remark.}
Equivalently, we could write the product of pre-free random variables
in the form
$$
[x_{m}^{(1)},q_{m}][x_{m}^{(2)},q_{m}]\ldots [x_{m}^{(n)},q_{m}]
=[x_{m+n-1}^{(1)}x_{m+n-2}^{(2)}\ldots x_{m}^{(n)},q_{m}],
$$
which is not hard to demonstrate. Heuristcally, the ranges of products of
pre-free random variables form a monotone increasing sequence. Note that
the same feature is exhibited by free random variables
acting in the free Fock space.\\[10pt]
\myownsection
\begin{center}
{\sc 4. Quantum semigroup structure on ${\cal F}({\cal G})$}
\end{center}
In this Section we show that one can endow
${\cal F}({\cal G})$ with a quantum semigroup (*-bialgebra) structure
with the algebraic tensor product $\otimes$ replaced by an appropriate `closure'
$\overline{\otimes}$.

This is done by first introducing a *-bialgebra structure on ${\cal F}_{0}({\cal G})$
and then lifting the comultiplication $\Delta$ and the counit $\epsilon$
from ${\cal F}_{0}({\cal G})$ to ${\cal F}({\cal G})$.
This can be done once a new type of tensor product, called
{\it monotone tensor product},
$$
{\cal F}({\cal G})\overline{\otimes}{\cal F}({\cal G})
=\overline{{\cal F}_{0}({\cal G})\otimes
{\cal F}_{0}({\cal G})}
$$
is introduced, where the monotone closure is taken on
strongly dense domains $(r_{m})$ implemented by the product lattice
${\cal P}^{(2)}=\Delta({\cal P})$

First, let us introduce a *-bialgebra
structure on ${\cal F}_{0}({\cal G})$.\\
\indent{\par}
{\sc Proposition 4.1.}
{\it The unital *-algebra ${\cal F}_{0}({\cal G})$ becomes a *-bialgebra, when
equipped with the comultiplication $\Delta:{\cal F}_{0}({\cal G})\rightarrow
{\cal F}_{0}({\cal G})\otimes {\cal F}_{0}({\cal G})$ given by
\begin{eqnarray*}
\Delta (X'(k))&=&X'(k)\otimes q_{k}+q_{k}\otimes X'(k)\\
\Delta (X''(k))&=&X''(k)\otimes q_{k-1}+q_{k-1}\otimes X''(k)\\
\Delta (q_{k})&=&q_{k}\otimes q_{k}\\
\Delta (1) &=& 1 \otimes 1
\end{eqnarray*}
and the counit $\epsilon:{\cal F}_{0}({\cal G}) \rightarrow \mathbb{C}$ given by
$\epsilon (X'(k))=\epsilon (X''(k))=0$
and $\epsilon (q_{k})=\epsilon (1)=1$, where $X\in {\cal G}$.}\\[5pt]
{\it Proof.}
It is easy to verify that $\Delta$ is coassociative
on generators since $X'(k)$ is $q_{k}$-primitive, $X''(k)$ is
$q_{k-1}$-primitive and $q_{k}$ as well as the unit $1$ are group-like.
Then the coassociativity on all of ${\cal F}_{0}({\cal G})$ easily follows.
Moreover, $\Delta$ and $\epsilon $ preserve the relations
in ${\cal F}_{0}({\cal G})$. For instance, if $k<m$, then
\begin{eqnarray*}
\Delta (q_{m}X'(k))&=&(q_{m}\otimes q_{m})(X'(k)\otimes q_{k}+
q_{k}\otimes X'(k))\\
&=&q_{m}X'(k)\otimes q_{k} + q_{k}\otimes q_{m}X'(k)\\
&=&X'(k)\otimes q_{k} + q_{k}\otimes X'(k)\\
&=&\Delta (X'(k))
\end{eqnarray*}
thus (3.2) is preserved (an identical proof holds for (3.3)).
In turn, if $k>m$, then
\begin{eqnarray*}
\Delta (q_{m}(X'(k)-X''(k)) &=&
q_{m}X'(k)\otimes q_{m}+ q_{m}\otimes q_{m}X'(k)\\
&-&
q_{m}X''(k)\otimes q_{m} - q_{m}\otimes q_{m}X''(k)\\
&=&
0
\end{eqnarray*}
an thus (3.4) is preserved.
Besides, it is easy to show that $\epsilon$ is a counit.
Thus, the triple $({\cal F}_{0}({\cal G}),\Delta, \epsilon)$ is a unital
*-bialgebra.\hfill $\Box$\\
\indent{\par}
In order to lift this structure to the monotone closure
${\cal F}({\cal G})$, we need to define a new type of tensor product called the
`monotone tensor product'. Our definition follows the pattern of the von Neumann algebra tensor product.

Let $({\cal B}^{(m)})$ and $({\cal C}^{(m)})$ be
sequences of increasing unital *-algebras (2.1) with
unions ${\cal B}_{0}$ and ${\cal C}_{0}$ (2.2), respectively
and let
$$
{\cal P}=\{b_{m};\; 0\leq m\leq \infty\},\;\;
{\cal Q}=\{c_{m};\; 0\leq m\leq \infty\},
$$
where $b_{0}=c_{0}=0$ and $b_{\infty}=1_{{\cal B}_{0}}$,
$c_{\infty}=1_{{\cal C}_{0}}$,
be the associated totally ordered lattices (2.4)
of projections which generate ${\cal B}^{(0)}$ and
${\cal C}^{(0)}$, respectively.
The sequence $({\cal B}^{(m)}\otimes {\cal C}^{(m)})$
is then an increasing sequence of unital *-algebras, for which we can construct
the monotone closure.
Thus, a MSDD in ${\cal B}_{0}\otimes {\cal C}_{0}$ is
a sequence of projections $(r_{m})$ from the lattice
$$
{\cal L}^{(2)}=\{b_{m}\otimes c_{m};\; 0\leq m\leq \infty\}
$$
such that $r_{m}\uparrow 1\otimes 1$.
In turn, an OWMC associated with ${\cal B}_{0}\otimes {\cal C}_{0}$ is
a sequence $(z_{m},r_{m})$, where $(r_{m})$ is a MSDD
in ${\cal B}_{0}\otimes {\cal C}_{0}$ and
$$
z_{m}=\sum_{k}w_{m}^{(k)}\otimes u_{m}^{(k)}\in {\cal B}_{0}^{(m)}
\otimes {\cal C}_{0}^{(m)}, \;\; m\geq 1
$$
for every $m$ (sums over $k$ are finite) and such that
$z_{n}r_{m}=z_{m}r_{m}$ and $z_{n}^{*}r_{m}=z_{m}^{*}r_{m}$
for all $n>m$.
An equivalence relation in the set of all such OWMC is analogous to that
of OWMC associated with ${\cal B}_{0}$.\\
\indent{\par}
{\sc Definition 4.1.}
A MCO `affiliated' with ${\cal B}_{0}\otimes {\cal C}_{0}$
is the equivalence class $[z_{m},r_{m}]$ corresponding
to OWMC $(z_{m},r_{m})$.
We denote by $\overline{{\cal B}_{0}\otimes {\cal C}_{0}}$
the unital *-algebra of MCO `affiliated'
with ${\cal B}_{0}\otimes{\cal C}_{0}$.
This closure leads to the definition of a
{\it monotone tensor product} denoted
${\cal B}\overline{\otimes} {\cal C}
=\overline{{\cal B}_{0} \otimes  {\cal C}_{0}}$.\\
\indent{\par}
{\it Remark 1.}
Note that by setting
$$
[x_{m},e_{m}]\otimes [y_{m},f_{m}]:=[x_{m}\otimes y_{m},e_{m}\otimes f_{m}]
$$
we obtain a natural *-algebra embedding
of ${\cal B}\otimes{\cal C}$ in ${\cal B}\overline{\otimes}{\cal C}$.\\
\indent{\par}
{\it Remark 2.}
For instance, the product in ${\cal B}\overline{\otimes}{\cal C}$
is given by
$$
[z_{m},r_{m}][w_{m}, s_{m}]
=[z_{m}w_{m}, k_{m}]
$$
where $k_{m}=r_{m}s_{m}w_{m}^{-1}(r_{m})(z_{m}^{*})^{-1}(s_{m})$
and $w_{m}^{-1}(r_{m})$ is defined to be the largest projection
$g\in {\cal L}^{(2)}$ for which $r_{m}w_{m}g=w_{m}g$
and thus, it always exists (similarly,
$(z_{m}^{*})^{-1}(r_{m})$ exists). \\
\indent{\par}
{\it Remark 3.}
One can proceed in a similar manner with
monotone tensor products of higher order
and this procedure is associative.
In particular
$$
{\cal B}\overline{\otimes}({\cal C}\overline{\otimes}{\cal D})=
({\cal B}\overline{\otimes} {\cal C})\overline {\otimes}{\cal D}=
{\cal B}\overline{\otimes}{\cal C}\overline{\otimes}{\cal D}
=
\overline{{\cal B}_{0}\otimes {\cal C}_{0}\otimes {\cal D}_{0}},
$$
with the closure taken w.r.t. MSDD from the lattice
$$
{\cal L}^{(3)}=
\{b_{m}\otimes c_{m}\otimes d_{m}; 0\leq m\leq \infty\}
$$
where the lattice associated with ${\cal D}_{0}$ is generated
by an increasing sequence of projections $(d_{m})$, with $d_{0}=0$ and
$d_{\infty}=1$.\\
\indent{\par}
Let us specify now how we can lift unital *-homomorphisms from
unital *-algebras to their monotone closures. Namely, let
${\cal B}_{0}$ and ${\cal C}_{0}$ be unital *-algebras
for which the associated monotone closures are
${\cal B}$ and ${\cal C}$, respectively.
If $\tau: {\cal B}_{0}\rightarrow {\cal C}_{0}$ is a unital
*-homomorphism, then one can lift $\tau $ to a unital
*-homomorphism from ${\cal B}$ to ${\cal C}$ by the formula
$$
\tau[x_{m},e_{m}]=[\tau (x_{m}),\tau (e_{m})],
$$
where we adopt the convention that the extended map is denoted
by the same symbol --
here, we understand that the monotone closure in ${\cal C}$
is `compatible' with $\tau$, i.e. is
taken w.r.t. MSDD from the lattice $\tau({\cal P})$.

In particular, if $\Delta $ is a comultiplication and $\epsilon$
is a counit on ${\cal B}_{0}$, then we can lift these maps to
${\cal B}$ (in the case of the counit, $\epsilon ({\cal P})=\{0,1\}$).
In order to treat the axioms required from the comultiplication
and the counit, we need to define a monotone tensor
product of maps
$\tau:{\cal B}_{0}\rightarrow {\cal D}_{0}$ and
$\sigma:{\cal C}_{0}\rightarrow {\cal E}_{0}$, namely
$$
\tau \overline{\otimes} \sigma : {\cal B}\overline{\otimes} {\cal C}
\rightarrow {\cal D} \overline{\otimes} {\cal E}
$$
$$
{\tau} \overline{\otimes} {\sigma}
[\sum_{k}w_{m}^{(k)}\otimes u_{m}^{(k)}, e_{m}\otimes f_{m}]
=
[\sum_{k}\tau(w_{m}^{(k)})\otimes \sigma (u_{m}^{(k)}),
\tau(e_{m})\otimes \sigma (f_{m})]
$$
where the RHS is a MCO in
${\cal D}\overline{\otimes} {\cal E}$. In this way we define
${\Delta} \overline{\otimes} {\rm id}$ and
${\rm id} \overline{\otimes} {\Delta}$ and their compositons with
$\Delta$, $({\Delta} \overline{\otimes} {\rm id})\circ \Delta$ and
$({\rm id} \overline{\otimes} {\Delta})\circ \Delta$ needed for
coassociativity.

Using these preparations, we can lift comultiplications and counits from
unital *-algebras to their monotone closures. The same axioms
hold as in the algebraic case except that the algebraic tensor
product $\otimes $ is replaced by the monotone tensor product
$\overline{\otimes}$. Thus we obtain `a unital *-bialgebra with respect to
the monotone tensor product' which we call a {\it monotone closed quantum semigroup}.

Let us now look at the case of ${\cal F}({\cal G})$. Note that since
all elements of ${\cal P}$ are in this case group-like w.r.t.
$\Delta$ and ${\cal P}^{(m)}=\Delta^{(m-1)} ({\cal P})$,
we can view ${\cal P}$ as a `group-like lattice'. Let us now define
for each natural number $n$ the lattice of projections
$$
{\cal P}^{(n)}=
\Delta^{(n-1)} ({\cal P}):
=\{\Delta^{(n-1)}(p), p\in {\cal P}\}
$$
where $\Delta^{(n)}$ is the $n$-th iteration
of the comultiplication $\Delta$, namely
$$
\Delta^{(n)}:=({\rm id}\otimes \Delta^{(n-1)})\circ \Delta,
$$
which allows us to compare projections from ${\cal P}^{(n)}$
for each $n$ with the order inherited from ${\cal P}$, i.e.
if $r=\Delta^{(n-1)}(p), s=\Delta^{(n-1)}(q)$, then
$r<s $ iff $p<q$.\\
\indent{\par}
{\sc Theorem 4.2.}
{\it The unital *-algebra ${\cal F}({\cal G})$ becomes a monotone closed
quantum semigroup when equipped with the comultiplication
$$
{\Delta}:
{\cal F}({\cal G}) \rightarrow {\cal F}({\cal G})
\overline{\otimes}
{\cal F}({\cal G}), \;\;
$$
$$
{\Delta}[x_{m},q_{m}]
:=
[\Delta (x_{m}), \Delta (q_m)],
$$
and the counit
${\epsilon}:{\cal F}({\cal G})\rightarrow \mathbb{C} $ given by
${\epsilon}([x_{m},e_{m}]):=[\epsilon (x_{m}),\epsilon (e_{m})]$.}\\[5pt]
{\it Proof.}
Coassociativity of $\Delta$ on
${\cal F}({\cal G})$
follows from the coassociativity of $\Delta$ on ${\cal F}_{0}({\cal G})$.
Verification of this fact and of the axioms on the counit is routine.
\hfill $\Box$\\
\indent{\par}
We will return to the monotone closed semigroup structure on ${\cal F}({\cal G})$ in Section 6. By then it will
become clear why this structure can be used to implement the free additive convolution. In the meantime, we will
study the tensor product representation of the free product of states in the general case of arbitrary unital
*-algebras.\\[10pt]
\myownsection
\begin{center}
{\sc 5. Free products}
\end{center}
In this section we will concentrate on reconstructing the free product of states from
the tensor product of states.

In order to do this we replace the `rigid' requirement for $\Delta$ to be a homomorphism by the requirement that
the `half-coproduct' maps of type (1.1)-(1.3) be homomorphisms.
For simplicity, we first establish our result for
two unital *-algebras and then generalize it to the case of an arbitrary family of unital
*-algebras.

Let ${\cal A}$ and ${\cal B} $ be arbitrary unital *-algebras.
Let us take free products
$$
\widehat{\cal H}({\cal A})=*_{k\in {\bf N}}{\cal A}(k),\;\; \widehat{\cal H}({\cal B})=*_{k\in {\bf N}}{\cal B}(k)
$$
of copies of ${\cal A}$ and ${\cal B}$ (in these free products we identify units), respectively,
and extend them by sequences of orthogonal projections
\begin{eqnarray*}
\widehat{\cal H}_{0}({\cal A})&=&\widehat{\cal H}({\cal A})*{\mathbb C}[p_{1},p_{2},p_{3},\ldots ]\\
\widehat{\cal H}_{0}({\cal B})&=&\widehat{\cal H}({\cal B})*{\mathbb C}[p_{1}',p_{2}',p_{3}',\ldots ]
\end{eqnarray*}
with the associated lattices of increasing projections
$$
{\cal P}_{1}=\{q_{m}: 0\leq m \leq \infty\}, \;\;
{\cal P}_{2}=\{q_{m}': 0\leq m \leq \infty\}
$$
obtained from the sequences $(p_{k})$ and $(p_{k}')$ by addition as in (2.3).

Next, let $I_{1}$ and $I_{2}$ be the two-sided *-ideals
in $\widehat{\cal H}_{0}({\cal A})$ and $\widehat{\cal H}_{0}({\cal B})$ generated by
$(1-q_{m})X(k)$ and $(1-q_{m}')Y(k)$ for $k<m$, respectively,
where $X(k)$ and $Y(k)$ denote the $k$-th copies of $X\in {\cal A}$ and $Y\in {\cal B}$.
Denote by ${\cal H}_{0}({\cal A})=\widehat{\cal H}_{0}({\cal A})/I_{1}$, ${\cal H}_{0}({\cal B})=
\widehat{\cal H}_{0}({\cal B})/I_{2}$
the corresponding quotient algebras and by ${\cal H}({\cal A})$, ${\cal H}({\cal B})$ their monotone closures.

Let us introduce the linear mappings
$$
j_{1}:{\cal A}\rightarrow {\cal H}({\cal A})\overline{\otimes}{\cal H}({\cal B}), \;\;
j_{2}:{\cal B}\rightarrow {\cal H}({\cal A})\overline{\otimes}{\cal H}({\cal B})
$$
given by
\begin{eqnarray}
j_{1}(X)&=&\sum_{k=1}^{\infty}X(k)\overline{\otimes} p_{k}'\\
j_{2}(Y)&=&\sum_{k=1}^{\infty}p_{k}\overline{\otimes} Y(k)
\end{eqnarray}
where $X \in {\cal A}$ and $Y \in {\cal B}$, where the associated MSDD is $(q_{m}\otimes q_{m}')$.
From now on we adopt the convention that monotone closed operators from the tensor product
${\cal H}({\cal A})\overline{\otimes}{\cal H}({\cal B})$, when written in the form of infinite sums, are
taken with MSDD of this form. \\
\indent{\par}
{\sc Proposition 5.1.}
{\it The mappings $j_{1}$ and $j_{2}$ are unital *-homomorphisms}.\\
\indent{\par}
{\it Proof.}
The mappings $j_{1}$, $j_{2}$ are *-homomorphisms due to orthogonality of the
sequences $(p_{k})$, $(p_{k}')$. Moreover, they are unital since, for instance
$$
j_{1}(1)=
[\sum_{k=1}^{m}
1(k)\otimes p_{k}', q_{m}\otimes q_{m}']= [1\otimes q_{m}', q_{m}\otimes q_{m}']=[1\otimes 1', 1\otimes 1']\equiv {\bf 1}
$$
where $1$ and $1'$ are units in ${\cal H}({\cal A})$ and ${\cal H}({\cal B})$, respectively.
A similar proof holds for $j_{2}$.\hfill $\Box$\\
\indent{\par}
Heuristically, the lemma below says that a certain `conditional expectation'
of alternating mixed moments of $j_{1}(X)$'s and $j_{2}(X)$'s takes values
in the two-sided ideal generated by `singletons'. We  shall use
\begin{eqnarray*}
{\cal K}&=&\langle q_{k}X(k)q_{k},\; {\rm where}\;k\geq 1\; {\rm and}\;X\in {\cal A}\rangle\\
{\cal K}'&=&\langle q_{k}'Y(k)q_{k}',\; {\rm where}\;k\geq 1\; {\rm and}\;Y\in {\cal B}\rangle ,
\end{eqnarray*}
the monotone closed two sided ideals in ${\cal H}({\cal A})$ and ${\cal H}({\cal B})$, respectively.\\
\indent{\par}
{\sc Lemma 5.2.}
{\it Let $E: {\cal H}({\cal A})\overline{\otimes}{\cal H}({\cal B})\rightarrow {\cal H}({\cal A})\overline{\otimes}
{\cal H}({\cal B})$ be given by}
$$
E(w)= (q_{1}\overline{\otimes} q_{1}')w(q_{1}\overline{\otimes} q_{1}')
$$
{\it and let $X_{k}\in {\cal A}_{i_{k}}$, where $k=1, \ldots , n$ and
$i_{1}\neq i_{2}\neq \ldots\neq i_{n}$ with ${\cal A}_{1}={\cal A}$, ${\cal A}_{2}={\cal B}$. Then}
$$
E(j_{i_{1}}(X_{1})j_{i_{2}}(X_{2})\ldots j_{i_{n}}(X_{n})) \subset {\cal I}
$$
{\it where ${\cal I}={\cal K}\overline{\otimes} {\cal H}({\cal B})+
{\cal H}({\cal A})\overline{\otimes} {\cal K}'$}.\\[5pt]
{\it Proof.}
To fix attention, we suppose that $i_{n}=1$ and denote
$$
A_{n}=j_{i_{1}}(X_{1})j_{i_{2}}(X_{2})\ldots j_{i_{n}}(X_{n})).
$$
We claim that
\begin{equation}
A_{n}(q_{1}\overline{\otimes} q_{1})=a_{n}\overline{\otimes} a_{n}'
\end{equation}
where, if $n$ is odd,
\begin{eqnarray*}
a_{n}&=&X_{1}(n)X_{3}(n-2)\ldots X_{n}(1)q_{1}\\
a_{n}'&=&X_{2}(n-1)X_{4}(n-3)\ldots X_{n-1}(2)q_{1}'
\end{eqnarray*}
and, if $n$ is even,
\begin{eqnarray*}
a_{n}&=&X_{2}(n-1)X_{4}(n-3)\ldots X_{n}(1)q_{1}\\
a_{n}'&=&X_{1}(n)X_{3}(n-2)\ldots X_{n-1}(2)q_{1}'
\end{eqnarray*}
for all $X_{n}, X_{n-2}, \ldots \in {\cal A}$ and $X_{n-1}, X_{n-3}, \ldots \in {\cal B}$.
We will use induction to prove the claim.
For $n=1$, the assertion is true since in that case
$A_{1}(q_{1}\overline{\otimes} q_{1}')= X_{1}(1)q_{1}\overline{\otimes} q_{1}'$.
We suppose now that $n$ is odd and that
the claim holds for the product
of $n-1$ factors. Then
$$
A_{n}(q_{1}\overline{\otimes} q_{1}')
=
\sum_{k=1}^{\infty}(X_{1}(k)\overline{\otimes} p_{k}')(w_{1}\overline{\otimes} w_{2})
$$
where
\begin{eqnarray*}
w_{1}
&=&
X_{3}(n-2)X_{5}(n-4)\ldots X_{n}(1)q_{1}\\
w_{2}
&=&
X_{2}(n-1)X_{4}(n-3)\ldots X_{n-1}(2)q_{1}'.
\end{eqnarray*}
Let us analyze the contribution from every term indexed by $k$, namely
$$
B_{k}=u_{k}\overline{\otimes} (v_{k}-v_{k}')
$$
where
\begin{eqnarray*}
u_{k}
&=&
X_{1}(k)X_{3}(n-2)\ldots X_{n}(1)q_{1}\\
v_{k}
&=&
q_{k}'X_{2}(n-1)X_{4}(n-3)\ldots X_{n-1}(2)q_{1}'\\
v_{k}'
&=&
q_{k-1}'X_{2}(n-1)X_{4}(n-3)\ldots X_{n-1}(2)q_{1}'.
\end{eqnarray*}
\indent{\par}
{\it Case 1.}
If $k>n$, then we apply (3.2) and (3.3) to the effect that
\begin{eqnarray*}
q_{k}'X_{2}(n-1)
&=&
X_{2}(n-1)\\
q_{k-1}'X_{2}(n-1)
&=&
X_{2}(n-1)
\end{eqnarray*}
which gives $v_{k}'=v_{k}'$ and thus $B_{k}=0$.
\indent{\par}
{\it Case 2.}
In turn, if $k=n$, then
\begin{eqnarray*}
v_{n}&=& X_{2}(n-1)X_{4}(n-3)\ldots X_{n-1}(2)q_{1}'\\
v_{n}'&=& q_{n-1}X_{2}(n-1)q_{n-1}X_{4}(n-3)\ldots X_{n-1}(2)q_{1}'
\end{eqnarray*}
since
\begin{eqnarray*}
q_{n}'X_{2}(n-1)&=&X_{2}(n-1)\\
q_{n-1}'X_{4}(n-3)&=&X_{4}(n-3)
\end{eqnarray*}
respectively, which implies that
$$
u_{n}\overline{\otimes} v_{n}'\in {\cal J}
$$
and $u_{n}\overline{\otimes}v_{n}$ is of the form given in the claim.
\indent{\par}
{\it Case 3.}
Finally, if $k<n$, then
\begin{eqnarray*}
u_{k}\overline{\otimes} v_{k}
&=&
u_{k}\overline{\otimes} q_{k}'q_{n-1}'X_{2}(n-1)q_{n-1}'w\in {\cal J}\\
u_{k}\overline{\otimes} v_{k}'
&=&
u_{k}\overline{\otimes} q_{k-1}'q_{n-1}'X_{2}(n-1)q_{n-1}'w'\in {\cal J}
\end{eqnarray*}
for some words $w,w'$ which start with $X_{4}(n-3)$, using similar arguments as above.
This finishes the proof of the claim.

After multiplying (5.3) from the left by $q_{1}\overline{\otimes} q_{1}'$, we can see that
$$
E(A_{n})\subset {\cal J}
$$
since in the expression
$$
q_{1}X_{1}(n)X_{3}(n-2)\ldots X_{n}(1)q_{1},
$$
we can write $q_{1}=q_{1}q_{n}$ and $X_{3}(n-2)=q_{n}X_{3}(n-2)$ to produce
$q_{n}X_{1}(n)q_{n}$.
This completes the proof.
\hfill $\Box$\\
\indent{\par}
Lemma 5.2 paves the way for the reconstruction of the free product of states. It now suffices to define
states $\widehat{\mu}$ and $\widehat{\nu}$ on ${\cal H}({\cal A})$ and ${\cal H}({\cal B})$, respectively,
in such a way that the `algebraic singletons' of Lemma 5.2 become `probabilistic singletons'
w.r.t. $\widehat{\mu}\overline{\otimes}\widehat{\nu}$.
This leads to our main result, namely that the free product of states $\mu * \nu$ on the free product ${\cal A}*{\cal B}$ of unital
*-algebras can be represented by the monotone tensor product $\widehat{\mu}\overline{\otimes}\widehat{\nu}$
of extended states on ${\cal H}({\cal A})\overline{\otimes}{\cal H}({\cal B})$. For this purpose we use
the unital *-homomorphism $j=j_{1}*j_{2}$ from the free product ${\cal A}*{\cal B}$ with identified units
into the monotone tensor product of ${\cal H}({\cal A})$ and ${\cal H}({\cal B})$, where $j_{1}$ and $j_{2}$
are given by (5.1)-(5.2).

Roughly speaking, the states $\widehat{\mu}$ and $\widehat{\nu}$ will turn out to be tensor products of
boolean extensions of $\mu$ and $\nu$, respectively. Let us recall that the {\it boolean extension} of
a state $\mu$ on a unital *-algebra ${\cal A}$ to the unital *-algebra
$$
\widetilde{\cal A}= {\cal A}*{\mathbb C}[P]
$$
where $P$ is a projection, is a state $\widetilde{\mu}$ such that
$$
\widetilde{\mu}(wPv)=\widetilde{\mu}(w)\widetilde{\mu}(v),\;\; \widetilde{\mu}(P)=1
$$
for any $w,v\in \widetilde{\cal A}$, and which, when restricted to ${\cal A}$, agrees with $\mu$.\\
\indent{\par}
{\sc Theorem 5.3.}
{\it For given states $\mu$ and $\nu$ on ${\cal A}$ and ${\cal B}$, respectively, there exist
states $\widehat{\mu}$ and $\widehat{\nu}$ on ${\cal H}({\cal A})$ and ${\cal H}({\cal B})$
and a unital *-homomorphism}
$$
j:\;{\cal A}*{\cal B}\rightarrow {\cal H}({\cal A})\overline{\otimes} {\cal H}({\cal B})
$$
{\it such that}
$$
\mu * \nu = (\widehat{\mu}\overline{\otimes} \widehat{\nu})\circ j\;.
$$
{\it Thus, the unital *-algebras $j({\cal A})$ and $j({\cal B})$ are free with respect to
$\widehat{\mu}\overline{\otimes} \widehat{\nu}$.}\\
\indent{\par}
{\it Proof.}
We divide the proof into 3 steps.

{\it Step 1.}
First, we construct suitable states on ${\cal H}_{0}({\cal A})$ and
${\cal H}_{0}({\cal B})$ associated with states $\mu$ on ${\cal A}$ and
$\nu$ on ${\cal B}$. Extend ${\cal A}$ and ${\cal B}$ by projections $P$ and $P'$, respectively, namely
$$
\widetilde{\cal A}={\cal A}* \mathbb{C}[P], \;\;\widetilde{\cal B}={\cal B}*\mathbb{C}[P']
$$
and extend $\mu$ and $\nu$ to their boolean extensions $\widetilde{\mu}$ and
$\widetilde{\nu}$. Now, take tensor products
$$
\widehat{\cal A}=\widetilde{\cal A}^{\otimes \infty},\;\;
\widehat{\cal B}=\widetilde{\cal B}^{\otimes \infty},\;\;
\widehat{\Phi}=\widetilde{\mu}^{\otimes \infty},\;\;
\widehat{\Psi}=\widetilde{\nu}^{\otimes \infty}
$$
and introduce unital *-homomorphisms
$$
\xi :\;{\cal H}_{0}({\cal A})\rightarrow \widehat{\cal A}, \;\; \eta:\;{\cal H}_{0}({\cal B})\rightarrow \widehat{\cal B}
$$
given by
\begin{eqnarray}
\xi (X(k))&=&1^{\otimes (k-1)}\otimes X \otimes 1^{\otimes \infty}\\
\xi (q_{m})&=&1^{\otimes (m-1)}\otimes P^{\otimes \infty}\\
\xi (1)&=& 1^{\otimes \infty}.
\end{eqnarray}
with analogous formulas for $\eta$. Then the linear functionals
$$
\widehat{\mu}=\Phi \circ \xi , \;\; \widehat{\nu}=\Psi \circ \eta
$$
are states on ${\cal H}_{0}({\cal A})$ and ${\cal H}_{0}({\cal B})$, respectively.

{\it Step 2.}
Let us show now that these states can be extended to states on ${\cal H}({\cal A})$ and ${\cal H}({\cal B})$,
respectively, and their tensor product to ${\cal H}({\cal A})\overline{\otimes}{\cal H}({\cal B})$ by taking pointwise limits.
Denoting them also by $\widehat{\mu}$ and $\widehat{\nu}$, we set
\begin{eqnarray}
\widehat{\mu}({\bf z})
&=&\lim_{m\rightarrow \infty}\widehat{\mu}(z_{m})\\
\widehat{\nu}({\bf u})
&=&\lim_{m\rightarrow \infty}
\widehat{\nu}(z_{m})\\
(\widehat{\mu}\overline{\otimes}\widehat{\nu})
({\bf w})
&=&
\lim_{m\rightarrow \infty}(\widehat{\mu}\otimes \widehat{\nu})
(w_{m})
\end{eqnarray}
where ${\bf z}=[z_{m},e_{m}]\in {\cal H}({\cal A})$, ${\bf u}=[u_{m},f_{m}] \in {\cal H}({\cal B})$
and ${\bf w}=[w_{m},r_{m}]\in {\cal H}({\cal A})\overline{\otimes}
{\cal H}({\cal B})$. In fact, let us show that the linear functionals given by (5.7)-(5.9)
are well-defined and that they are states.
We shall establish existence of the limits on the RHS of (5.7) only since the proofs for (5.8)-(5.9)
are similar.
We have $e_{m}=q_{k(m)}$ with $k(m)\uparrow \infty$.
Since $[z_{m},q_{k(m)}]$ is MCO, we have
$z_{n}q_{k(m)}=z_{m}q_{k(m)}$ for $n>m$.
Since there exists natural $s$ such that $k(m)>0$ for $m>s$, hence
$$
z_{n}q_{1}=z_{n}q_{k(m)}q_{1}=z_{m}q_{k(m)}q_{1}=z_{m}q_{1}
$$
for $n>m>s$. Therefore, there exists $s$ such that
$$
q_{1}z_{n}q_{1}=q_{1}z_{m}q_{1}
$$
for $n>m>s$, which implies that
$\lim_{m\rightarrow \infty} \widehat{\mu}(z_{m})$
exists since
$$
\widehat{\mu} (z_{m})=
\widehat{\mu} (q_{1}z_{m}q_{1})=
\widehat{\mu} (q_{1}z_{n}q_{1})=
\widehat{\mu}(z_{n})
$$
for $n>m>s$, where we used $\widetilde{\mu}(PuP)=\widetilde{\mu}(u)$
and the definition of $\widehat{\mu}$.
Moreover, it does not depend on the choice of representatives
of MCO. In fact, choose a different representative for
${\bf z}$, say $(z_{m}', q_{l(m)})$
instead of $(z_{m}, q_{k(m)})$ with the equivalence $(z_{m},q_{k(m)})
\equiv (z_{m}',q_{l(m)})$ implemented by some MSDD $(e_{m})=(q_{r(m)})$, i.e.
$z_{m}e_{m}=z_{m}'e_{m}$ for all $m$. Again, there exists natural $s$
such that $r(m)>0$ for $m>s$. A similar argument as that given above
gives $z_{m}q_{1}=z_{m}'q_{1}$ for $m>s$,
which gives independence of the limit
from the choice of representatives.
Finally, positivity and normalization of the extended functional
follow from its definition.

{\it Step 3.}
Finally, let us reconstruct the free product of states.
Let
$$
j:{\cal A}*{\cal B}\rightarrow {\cal H}({\cal A})\overline{\otimes}{\cal H}({\cal B})
$$
be the unital *-homomorphism given by
$$
j(X)=j_{1}(X),\; j(Y)=j_{2}(Y)
$$
where $X\in {\cal A}$ and $Y\in {\cal B}$. Use Proposition 5.1 and the fact that
$j_{1}(1_{{\cal A}})=j_{2}(1_{{\cal B}})$ to conclude that $j$ is well-defined.
We will show that $\mu * \nu$ agrees with $(\widehat{\mu}\overline{\otimes}\widehat{\nu})\circ j$.
Note that if $X\in {\rm ker}\mu$, where $\mu$ is a state
on ${\cal A}$, then $PXP\in {\rm ker}\widetilde{\mu}$, where
$\widetilde{\mu}$ is the boolean extension of $\mu$. A similar
statement holds for $\nu$.
Thus, if ${\cal K}_{\mu}$ and $K_{\nu}$  denote the two-sided ideals
$$
K_{\mu}=\langle q_{k}X(k)q_{k},\; {\rm where }\;k\in {\bf N}\; {\rm and}\;
X\in {\cal A}\cap {\rm ker}\mu \rangle
$$
$$
K_{\nu}=\langle q_{k}'Y(k)q_{k}',\; {\rm where }\;k\in {\bf N}\; {\rm and}\;
Y\in {\cal B}\cap {\rm ker}\nu \rangle
$$
respectively,
then
$$
{\cal K}_{\mu}\overline{\otimes}{\cal H}({\cal B})
+{\cal H}({\cal A})\overline{\otimes}{\cal K}_{\nu}
\subset {\rm ker}(\widehat{\mu}\;\overline{\otimes }\;\widehat{\nu}) .
$$
It is now enough to remark that
$$
\widehat{\mu} \overline{\otimes} \widehat{\nu}
=
(\widehat{\mu}\overline{\otimes}\widehat{\nu}) \circ E
$$
since $\widetilde{\mu}(PuP)=\widetilde{\mu}(u)$ and $\widetilde{\nu}(P'wP')=\widetilde{\nu}(w)$.
Using this and Lemma 5.2, we get
$$
\widehat{\mu} \overline{\otimes} \widehat{\nu}
(j_{i_{1}}(X_{1})\ldots j_{i_{n}}(X_{n}))=0
$$
if $X_{1}, \ldots , X_{n}$ are in the kernels of
$\mu$ or $\nu$ depending on whether we have the (alternating) indices $i_{1}, \ldots , i_{n}$ equal to
$1$ or $2$, respectively. This completes the proof.
\hfill $\Box$\\
\indent{\par}
{\sc Corollary 5.4.}
{\it For given states $\mu$ and $\nu$ on ${\cal A}$ and ${\cal B}$, the random variables $j_{1}(X)$ and
$j_{2}(Y)$, where $X\in {\cal A}$ and $Y\in {\cal B}$, have the same distributions as $X$ and $Y$, respectively, and
are free with respect to the tensor product state $\widehat{\mu}\overline{\otimes}\widehat{\nu}$.}\\
\indent{\par}
{\it Proof.}
This is an immediate consequence of Theorem 5.3.
\hfill $\Box$\\
\indent{\par}
{\it Remark.}
Thus formulas (5.1)-(5.2) give a tensor product representation of free random variables. Treating states $\widehat{\mu}$
and $\widehat{\nu}$ as canonical extensions of $\mu$ and $\nu$ to ${\cal H}({\cal A})$ and ${\cal H}({\cal B})$,
we can view this representation as the one in which information about freeness is shifted from states to variables.\\
\indent{\par}
The above theorem can be generalized to arbitrary
families $({\cal A}_{l})_{l\in L}$ of unital *-algebras.
Let ${\cal H}_{0}({\cal A}_{l})$ be the unital *-algebra constructed from copies of ${\cal A}_{l}$ as
shown in Section 3 and let ${\cal H}({\cal A}_{l})$ be the associated unital *-algebra
of monotone closed operators. By $(q_{k,l})_{k\in {\mathbb N}}$ for each $l\in L$
we denote the associated increasing sequence of projections. Finally, let
$$
\overline{\bigotimes _{l\in L}}{\cal H}({\cal A}_{l})=\overline{\bigotimes _{l\in L}{\cal H}_{0}({\cal A}_{l})}
$$
where the monotone closure on the right-hand side is taken with respect to the product lattice
$$
{\cal L}=\{\otimes _{l\in L}q_{m,l}:\; 0\leq m \leq \infty \}
$$
where we set $q_{0,l}=0$ and $q_{\infty, l}=1_{l}$ for every $l\in L$.\\
\indent{\par}
{\sc Theorem 5.5.}
{\it For a given family of states $(\mu_{l})_{l\in L}$ on unital *-algebras $({\cal A}_{l})_{l\in L}$
there exist states $(\widehat{\mu})_{l\in L}$ on $({\cal H}({\cal A}_{l}))_{l\in L}$, respectively, and a unital
*-homomorphism}
$$
j:\ *_{l\in L}{\cal A}_{l}\rightarrow \overline{\bigotimes _{l\in L}}{\cal H}({\cal A}_{l})
$$
{\it such that}
$$
*_{l\in L}\mu_{l}=   (\overline{\otimes}_{l\in L}\widehat{\mu}_{l}) \circ j
$$
{\it and thus the unital *-algebras $(j({\cal A}_{l}))_{l\in L}$ are
free with respect to the monotone tensor product state $\overline{\otimes}_{l\in L}\widehat{\mu}_{l}$.}\\
\indent{\par}
{\it Proof.} For $X\in {\cal A}_{m}$ we define $j(X)$ by the formula
\begin{equation}
j(X)=\sum_{k=1}^{\infty}X(k)\overline{\otimes} p_{k}
\end{equation}
according to the decomposition
$$
\overline{\bigotimes _{l\in L}}{\cal H}({\cal A}_{l})={\cal H}({\cal A}_{m})
\overline{\otimes}\overline{\bigotimes _{l\neq m}}{\cal H}({\cal A}_{l})
$$
where $p_{k}=q_{k}-q_{k-1}$ and
$$
q_{k}=\overline{\bigotimes _{l\neq m}}q_{k,l}.
$$
Then we extend $j$ to the free product $*_{l\in L}{\cal A}_{l}$ by taking the linear and multiplicative extension
of (5.10) for all $m\in L$. The states $(\widehat{\mu}_{l})_{l\in L}$ are defined
as
$$
\widehat{\mu}_{l}= \Phi_{l} \circ \xi_{l}
$$
where
$$
\Phi_{l}=\widetilde{\mu}^{\otimes \infty}
$$
are states on the tensor product algebras
$$
\widehat{\cal A}_{l}=\widetilde{\cal A}_{l}^{\otimes \infty}
$$
where $\widetilde{\cal A}_{l}={\cal A}_{l}*\mathbb{C}[P_{l}]$ and $\xi_{l}$ is given by the
same formula as (5.4)-(5.6) with $P$ replaced by $P_{l}$.
The infinite tensor products are understood to be taken with respect to the
family $\{1_{l},P_{l}\}$ [F-L-S], by which we mean that only finitely many
sites are occupied by elements different from $1_{l}$'s and $P_{l}$'s.
The remaining part of the proof is similar to that in the case of two algebras.
\hfill $\Box$.\\[10pt]
\myownsection
\begin{center}
{\sc 6. Additive free convolution}
\end{center}
We return to the monotone closed semigroup structure on ${\cal F}({\cal G})$ in order to show that
it implements the additive free convolution.
In other words, the comultiplication $\Delta$ on ${\cal F}({\cal G})$
maps pre-free random variables onto the sum of random variables which are free with respect to
a tensor product state.

By applying the comultiplication $\Delta$
to the pre-free random variable associated with $X\in {\cal G}$ one obtains the sum
$$
\Delta (\sum_{k=1}^{\infty}\delta X(k))=J_{1}(X)+J_{2}(X)
$$
of random variables
\begin{eqnarray}
J_{1}(X)&=&\sum_{k=1}^{\infty}(X'(k)\overline{\otimes} q_{k}-X''(k)\overline{\otimes} q_{k-1})\\
J_{2}(X)&=&\sum_{k=1}^{\infty}(q_{k}\overline{\otimes} X'(k)-q_{k-1}\overline{\otimes} X''(k))
\end{eqnarray}
from ${\cal F}({\cal G}) \overline{\otimes} {\cal F}({\cal G})$. Although they differ from $j_{1}(X)$ and $j_{2}(X)$
given by (5.1)-(5.2), this difference is not relevant in the weak sense, namely they turn out to be free with respect to
the tensor product state $\widehat{\mu}\overline{\otimes}\widehat{\nu}$ for suitable $\widehat{\mu}$, $\widehat{\nu}$.

The states $\widehat{\mu}$ and $\widehat{\nu}$ will be obtained by lifting the corresponding states on
the quotient ${\cal H}({\cal G})$ of ${\cal F}({\cal G})$ modulo the two-sided ideal generated by
$X'(k)-X''(k)$, $X\in {\cal G}$, $k\in {\mathbb N}$ -- by abuse of notation we shall also denote these lifted states
by $\widehat{\mu}$ and $\widehat{\nu}$.  \\
\indent{\par}
{\sc Theorem 6.1.}
{\it For given states $\mu$ and $\nu$ on ${\cal A}$ there exist states $\widehat{\mu}$ and $\widehat{\nu}$
on ${\cal F}({\cal G})$ and a unital *-homomorphism}
$$
J:\;{\cal A}*{\cal A}\rightarrow {\cal F}({\cal G})\overline{\otimes}{\cal F}({\cal G})
$$
{\it such that $\mu * \nu$ agrees with  $(\widehat{\mu} \overline{\otimes} \widehat{\nu})\circ J$
on ${\cal A}*{\cal A}$.}\\
\indent{\par}
{\it Proof.}
The proof is essentially the same as that of Theorem 5.2. The states $\widehat{\mu}$ and $\widehat{\nu}$ are now
defined by
$$
\widehat{\mu}=\Phi \circ \xi \circ i , ;\;\ \widehat{\nu}=\Psi \circ \eta \circ i
$$
where $i:{\cal F}({\cal G})\rightarrow {\cal H}({\cal G})$ is the identification map given by
$i(X'(k))=i(X''(k))=X(k)$ for every $X\in {\cal G}$ and $i(q_{k})=q_{k}$. Then we set $J=J_{1}*J_{2}$, i.e.
$J$ is the homomorphic extension of
$$
J(X)=\left\{
\begin{array}{ll}
J_{1}(X) & {\rm if}\; X\in {\cal G}_{1}\\
J_{2}(X) & {\rm if}\; X\in {\cal G}_{2}
\end{array}\right.
$$
and $J(1_{{\cal A}})=1\otimes 1$, where ${\cal G}_{1}$ and ${\cal G}_{2}$ denote copies of ${\cal G}$.
The rest of the proof is analogous to that of Theorem 5.2.
\hfill $\Box$\\
\indent{\par}
{\sc Corollary 6.2.}
{\it Let $\mu \boxplus \nu$ be the free additive convolution of states $\mu$ and $\nu$ on ${\cal A}$ and let
$\widehat{\mu}\star \widehat{\nu}=(\widehat{\mu}\overline{\otimes}\widehat{\nu})\circ \Delta$ be the convolution
of states of Theorem 6.1. Then}
$$
\mu \boxplus \nu
=
(\widehat{\mu}\star \widehat{\nu})\circ \tau ,
$$
{\it where $\tau:{\cal A}\rightarrow {\cal F}({\cal G})$
is the unital *-homomorphism given by}
$$
\tau(X)=\sum_{k=1}^{\infty}\delta X(k), \;\; \tau(1_{{\cal A}} )=1
$$
{\it where $X\in {\cal G}$.}
\\[5pt]
{\it Proof.}
This is an immediate consequence of Theorem 6.1.\hfill $\Box$\\
\indent{\par}
Thus we can use a `quantum group' language to speak about the additive free convolution instead
of using the dual group language. In particular, let ${\cal A}=\mathbb{C}[X]$, ${\cal G}=\{X\}$, and
let $\mu , \nu$ be states on $\mathbb{C}[X]$ corresponding to measures $\mu ,\nu$ on ${\mathbb R}$.
Heuristically, the monotone closed quantum semigroup ${\cal F}({\cal G})$
can be interpreted as a quantum analog of ${\mathbb R}^{\infty}$ as it is the algebra of
polynomials in countably many noncommuting variables.
Therefore, the additive free convolution of classical measures $\mu \boxplus\nu$
can be viewed as a restriction of the convolution $\widehat{\mu}\star \widehat{\nu}$ of
`quantum measures' on `quantum ${\mathbb  R}^{\infty}$' to the `quantum free line'
${\cal F}_{{\rm pf}}({\cal G})$.
.\\[10pt]
\begin{center}
{\sc Acknowledgements}
\end{center}
I would like to thank Professor Luigi Accardi for his remarks and suggestions which were very helpful in the
preparation of the revised version of the mansucript.\\[10pt]
\begin{center}
{\sc References}
\end{center}
[Av] {\sc D.~Avitzour}, ``Free products of $C^{*}$- algebras'',
{\it Trans.~Amer.~Math.~Soc.} {\bf 271} (1982), 423-465.\\[3pt]
[B] {\sc M.~Bozejko}, ``Uniformly bounded representations of free groups'',
{\it J. Reine Angew. Math.} {\bf 377} (1987), 170-186.\\[3pt]
[Be1] {\sc S.~Berberian}, {\it Baer *-Rings}, Springer-Verlag, 1972.\\[3pt]
[Be2] {\sc S.~Berberian}, ``The regular ring of a finite Baer *-ring'',
{\it J. Algebra} {\bf 23} (1972), 35-65.\\[3pt]
[Be3] {\sc S.~Berberian}, ``The regular ring of a finite AW$^{*}$- algebra'',
{\it Ann. Math.} {\bf 65} (1957), 224-240.\\[3pt]
[F-L-S] {\sc U.~Franz, R.~Lenczewski, M.~Sch\"{u}rmann},
``The GNS construction for the hierarchy of freeness'', Preprint No.
9/98, Wroclaw University of Technology, 1998.\\[3pt]
[L1] {\sc R.~Lenczewski}, ``Unification of independence in
quantum probability'', {\it Inf.~Dim. Anal.~Quant.~Probab. Rel.~Topics}
{\bf 1} (1998), 383-405.\\[3pt]
[L2] {\sc R.~Lenczewski}, ``Filtered random variables, bialgebras and
convolutions'', {\it J. Math. Phys.}
{\bf 42} (2001), 5876-5903.\\[3pt]
[M] {\sc N.~Muraki}, ``Monotonic independence, monotonic central limit theorem and monotonic law of large numbers'',
{\it Inf.~Dim.~Anal.~Quant.~Probab.~Rel.~Topics.} {\bf 4} (2001), 39-58.\\[3pt]
[V1] {\sc D.~Voiculescu}, ``Symmetries of some reduced free product $C^{*}$-
algebras'', in {\it Operator Algebras and their Connections with Topology
and Ergodic Theory}, Lect. Notes Math. {\bf 1132},
556-588 (1985).\\[3pt]
[V2] {\sc D.~Voiculescu}, ``Addition of certain non-commuting random
variables'', {\it J. Funct. Anal.} {\bf 66} (1986), 323-346.\\[3pt]
[V3] {\sc D.~Voiculescu}, ``Dual algebraic structures on operator algebras
related to free products'', {\it J. Operator Theory} {\bf 17} (1987),
85- 98.
\end{document}